\documentclass[10pt,a4paper]{article}
\usepackage{amsfonts}
\usepackage{amssymb}
\usepackage{mathrsfs}
\usepackage{amsthm}

\usepackage{fancyhdr}

\newtheorem{defi}{Definition}[section]
\newtheorem{theo}[defi]{Theorem}
\newtheorem{prop}[defi]{Proposition}

\newtheorem{lemm}[defi]{Lemma}

\theoremstyle{definition}
\newtheorem*{rem}{Remark}
\newtheorem{nrem}{Remark}[section]

\newcommand{\be}{\begin{eqnarray*}}
\newcommand{\ee}{\end{eqnarray*}}
\newcommand{\beqa}{\begin{eqnarray}}
\newcommand{\eeqa}{\end{eqnarray}}
\newcommand{\ba}{\begin{array}}
\newcommand{\ea}{\end{array}}

\newtheorem*{exa}{Example}

\begin{document}
\title{Definite signature conformal holonomy: a complete classification}
\author{Stuart Armstrong}
\date{12 September 2006}
\maketitle

\begin{abstract}
This paper aims to classify the holonomy of the conformal Tractor connection, and relate these holonomies to the geometry of the underlying manifold. The conformally Einstein case is dealt with through the construction of metric cones, whose Riemmanian holonomy is the same as the Tractor holonomy of the underlying manifold. Direct calculations in the Ricci-flat case and an important decomposition theorem complete the classification for definitive signature.
\end{abstract}

\tableofcontents

\subsection*{Acknowledgements:}
The Author would like to acknowledge the support of an EPSRC research studentship. He would like to thank his supervisor, Prof.~Nigel Hitchin, for all the help provided, and Dr.~Thomas Leistner for several illuminating exchanges.

\section{Introduction}

Conformal geometry is perhaps the most natural extension of Riemannian geometry, and shares many of the same features with it. However, it was realised early on -- as far back as Cartan \cite{ECC} -- that one of the most mathematically rewarding ways of dealing with conformal geometry was not by analogy with Riemannian geometries, but by analogy with the other parabolic geometries, using the general \emph{Cartan connection} as a universal tool.

These parabolic geometries are a class of geometries that include, amongst others, projective, almost Grassmanian, almost quaternionic, and co-dimension one CR structures. The common point of these is that their `flat' model space is the Lie group quotient $G/P$, where $P$ is parabolic. Papers \cite{mepro1} and \cite{mepro2} by the same author deals with the projective case, while this paper treats the conformal one.

Many figures contributed to understanding parabolic geometries; T.Y.~Thomas \cite{OCG}, \cite{CT} developed key ideas for Cartan connection calculus, and Shiego Sasaki investigated the conformal case in 1943 \cite{Sasaki}, \cite{Sasaki2}, followed by N.~Tanaka \cite{OEPA} in 1979 and the major paper of T.N.~Bailey, M.G.~Eastwood and R.~Gover in 1994 \cite{TSB}.

Since then, there have been a series of papers by A.~\v Cap and R.~Gover \cite{TCPG}, \cite{TBIPG}, \cite{CCG}, \cite{STCAMC}, developing a lot of the techniques that will be used in the present paper. Previous papers had focused on seeing the Cartan connection for conformal geometry as a property of a principal bundle $\mathcal{P}$. More recently, the principal bundle is replaced by an associated vector bundle, the \emph{Tractor bundle} $\mathcal{T}$, and the Cartan connection by a connection form for $\mathcal{T}$, the \emph{Tractor connection} $\overrightarrow{\nabla}$. With these tools, calculations are considerably simplified.

The purpose of this paper is to analyse one of the invariants of the Tractor connection, the holonomy group. There is an invariant metric of signature $(n+1, 1)$ on $\mathcal{T}$, so this holonomy group must be a sub-group of $G = SO(n+1,1)$.

It is a well known fact that a parallel section of the Tractor bundle corresponds to the local existence of an Einstein metric in the conformal class of a manifold. Beyond this, little was known about reductions of holonomy.

In this paper, we shall classify all the possible local holonomy groups of $\overrightarrow{\nabla}$ acting \emph{reducibly} on $\mathcal{T}$. In doing so, they must conserve a Lorentzian metric of signature $(n+1,1)$. Then a paper by A.J.~Di Scala and C.~Olmos \cite{GHSHS} shows that we have the complete list: there exist no connected proper subgroups of $SO(n+1,1)$ acting irreducibly on $\mathbb{R}^{n+1,1}$.
\begin{prop}
There are no local holonomy algebras acting irreducibly on $\mathcal{T}$ apart from the full $\mathfrak{so}(n+1,1)$ algebra.
\end{prop}

A very recent paper by Felipe Leitner, \cite{NCKF}, proves the same results as in this paper; but his methods, involving normal Killing Spinors, are different from those described here.

The classification comes in two main pieces; if a bundle of rank other than $1$ or $n$ is preserved, the manifold decomposes analogously to the De Rham decomposition:
\begin{theo} Let $(M, [g])$ be a conformal, $n$-dimensional manifold, such that $\mathcal{T}_M$ has a holonomy preserved sub-bundle of rank $k$. Then there exists a metric $g \in [g]$ such that $(M, g)$ splits locally into the direct product of two Einstein manifolds $N_1$, $N_2$ of dimensions $l = k-1$ and $n-l$. The Einstein constants $a$ and $b$ of $N_1$, $N_2$ are related by $(n-l-1)a = (1-l)b$. Furthermore, there are cannoincal inclusions of the Tractor bundles of $N_1$ and $N_2$ into $\mathcal{T}_M$ and the Tractor holonomy group of $M$ is the direct product of those of $N_1$ and $N_2$.\end{theo}
That last statement requires a bit of explaining, since the Tractor bundles of $N_1$ and $N_2$ are of rank $l +2$ and $n-l+2$ respectively. However, since these are both Einstein manifolds, the effective rank of their Tractor bundles are $l+1$ and $n-l+1$, allowing the decomposition.

The converse is also true. This decomposition is a local result, and may become degenerate along some embedded sub-manifolds.

The second step is to list all the possible Tractor holonomies for a conformally Einstein manifold. Using a metric cone construction, related to the Ambient Metric of \cite{CI}, \cite{STCAMC} and \cite{AOT}, the following list is established:
\begin{theo}[Einstein Classification] The Tractor holonomy of $(M^n, [g])$, for $M^n$ conformal to an Einstein space of non-zero scalar curvature, is one of the following groups:
\begin{itemize} \item[-] $SO(n,1), n \geq 4$, \item[-] $SO(n+1), n \geq 4$, \item[-] $SU(m)$ for $2m=n+1, n \geq 4$, \item[-] $Sp(m)$ for $4m=n+1$, \item[-] $G_2$ for $n=6$, \item[-] $Spin(7)$ for $n=7$. \end{itemize}
Moreover, all these actually occur as holonomy groups.
\end{theo}

The Ricci-flat case must be treated differently; in fact, if $(M^n,g)$ is Ricci-flat and conformally indecomposable, and $G$ is the metric holonomy group of $\nabla^g$, then $(M^n, [g])$ has Tractor holonomy $G \rtimes \mathbb{R}^n$. Thus:

\begin{theo} The possible indecomposable Tractor holonomy groups for the conformal manifold $(M^n, [g])$, conformally Ricci-flat, are:
\begin{itemize}
\item[-] $SO(n) \rtimes \mathbb{R}^n, n \geq 4$,
\item[-] $SU(m) \rtimes \mathbb{R}^{2m}, m \geq 2$,
\item[-] $Sp(m) \rtimes \mathbb{R}^{4m}, m \geq 1$,
\item[-] $G_2 \rtimes \mathbb{R}^7$,
\item[-] $Spin(7) \rtimes \mathbb{R}^8$,
\end{itemize}
and all of these groups do occur.
\end{theo}

This paper begins with defining and laying out the groundwork for the conformal Tractor Bundle and connection. Furthermore, it will prove the equivalence of this (second order) point of view with the standard view of the conformal structure as an equivalence class of metric structures. Some standard results will then be presented, showing how an Einstein structure in the conformal class is equivalent to a parallel section of the Tractor bundle.

Section 4, the heart of the paper, introduces umbilicity, the conformal equivalent of totaly geodicity, and proves the decomposition theorem previously mentioned.

Section 5 then establishes the list for the Einstein spaces via the metric cone construction, with Section 6 complementing it using different methods to list the possible holonomies for conformally Ricci-flat manifolds.

A brief note on symmetric spaces follows, to illustrate the use of these methods; the paper ends with considerations of the differences that arise with indefinite signature.

This paper formed the beginning of the author's thesis \cite{methesis} and was inspired and supervised by Dr. Nigel Hitchin.

\begin{rem} In all the holonomy groups listed in this paper, the holonomy reduction corresponds to the existence of a particular metric in the conformal class. Hence we always have a canonical representative in the conformal class, whenever the holonomy reduces. \end{rem}

\section{Cartan Connection: Theory}

\subsection{The Cartan Connection}

With homogenous geometries, since Klein, one deals with spaces $X = G/P$, for $G$ a Lie group acting transitively and effectively on $X$ and $P$ a subgroup.

The Cartan connection is a curved version of the flat geometries. Given \emph{any} manifold $M$, it maps the tangent space $T_M$ locally to the Lie algebra quotient,
\be
(T_M)_x \cong \mathfrak{g}/\mathfrak{p},
\ee
for all $x$ in $M$.

We will follow the exposition used in \cite{TCPG}. In all of the following, we assume that $M$ is an $n$-dimensional manifold, with $\mathfrak{g}$ a semisimple Lie algebra and a subalgebra $\mathfrak{p} \subset \mathfrak{g}$ with $\mathfrak{p}$ of codimension $n$ in $\mathfrak{g}$. There are corresponding groups $P \subset G$; different choices of such groups may change the global properties of Cartan connections, but not the local ones.

\begin{defi}[Cartan Connection] \label{NCC:def} On $M$, given a principal $P$-bundle $\mathcal{P} \to M$, a normal Cartan connection $\omega$ is a section of $T^*_{\mathcal{P}} \otimes \mathfrak{g}$, with the following properties:
\begin{enumerate}
\item $\omega$ is invariant under the $P$-action ($P$ acting by $Ad$ on $\mathfrak{g}$),
\item $\omega(\sigma_A) = A$, where $\sigma_A$ is the fundamental vector field of $A \in \mathfrak{p}$,
\item $\omega_u: T \mathcal{P}_u \to \mathfrak{g}$ is a linear isomorphism for all $u \in \mathcal{P}$.
\end{enumerate}
\end{defi}
If $\mathfrak{p}$ is a moreover a parabolic subalgebra (see paper \cite{riccor}), we may make the further requirement that the connection be \emph{normal}; this is a uniqueness condition for the Cartan connection of a particular geometry, similar to the torsion-free condition for a Levi-Civita connection. See \cite{TCPG} for a proof of the existence of a normal Cartan connection in all parabolic geometries.

Paper \cite{riccor} defines a parabolic subalgebra in an elegant and invariant way; for our purposes, however, it suffices to require that there exist a graded splitting of $\mathfrak{g}$
\be
\mathfrak{g} = \mathfrak{g}_{l} \oplus \mathfrak{g}_{l-1} \oplus \ldots \oplus \mathfrak{g}_{1} \oplus  \mathfrak{g}_{0} \oplus  \mathfrak{g}_{-1} \oplus \ldots \oplus \mathfrak{g}_{-l},
\ee
such that $[\mathfrak{g}_{j}, \mathfrak{g}_{k} ] \subseteq \mathfrak{g}_{j+k}$ and 
\be
\mathfrak{p} = \mathfrak{g}_{0} \oplus \mathfrak{g}_{-1} \ldots \oplus \mathfrak{g}_{-l}.
\ee
The algebra is then called $|l|$-graded; the conformal algebra will be seen to be $|1|$-graded.

\begin{defi}[Normal Cartan Connection]
A Cartan connection for a given parabolic geometry is normal if it has the following additional condition:
\begin{enumerate}
\setcounter{enumi}{3}
\item The `curvature' $\kappa (\eta,\xi) = d \omega (\eta,\xi) + [\omega(\eta), \omega(\xi)]$ is such that $\partial^* \kappa = 0$ where $\partial^*$ is the dual homology operator.
\end{enumerate}
\end{defi}
There is, however, a simpler characterisation of the normality condition in the conformal case, see the proof of Lemma (\ref{norm:lemm}).

The bundle $\mathcal{P}$ and the form $\omega$ together define the geometry. The first two conditions on $\omega$ are analogous to those of a standard connection. The third condition is very different, however, giving a pointwise isomorphism $T\mathcal{P}_u \to \mathfrak{g}$ rather than a map with kernel.

However the Cartan connection does give rise to a connection in the usual sense, the so-called Tractor connection.

The inclusion $P \hookrightarrow G$ generates a principal bundle inclusion $i: \mathcal{P} \hookrightarrow \mathcal{G}$, with $\mathcal{G}$ a $G$-bundle, and generates a standard connection form:

\begin{prop} \label{exten:conn} There is a unique $\omega' \in \Omega^1(\mathcal{G}, \mathfrak{g})$ such that $\omega'$ is a standard connection form on $\mathcal{G}$ and $i^* \omega' = \omega$. \end{prop}
\begin{proof}

At any point of $\mathcal{P} \hookrightarrow \mathcal{G}$, define $\omega'(X) = \omega(X)$ for $X \in \Gamma(T\mathcal{P})$, and $\omega'(\sigma_A) = A$ for $\sigma_A$ the fundamental vector field of $A \in \mathfrak{g}$. These two formulas correspond whenever they are both defined (Property 2 from Definition \ref{NCC:def}), and completely define $\omega'$ on $\mathcal{P}$. Then define $\omega'_u = g^*(\omega'_{g(u)})$ in the general case, for $g(u) \in \mathcal{P}$. Property 1 for $\omega$ ensures this is well defined.

To see that $\omega'$ is indeed a connection, notice that for $v \in \mathcal{P}$, $\omega' : T\mathcal{G}_v \to \mathfrak{g}$ has maximal rank, since $\omega = \omega'|_{T\mathcal{P}} : T\mathcal{P}_v \to \mathfrak{g}$ is surjective. $G$-invariance of $\omega'$ generalises this property to all of $\mathcal{G}$. \end{proof}

This $\omega'$ is the Tractor connection; when we see it as a connection on an associated vector bundle, we shall designate it by $\overrightarrow{\nabla}$. The Tractor connection obviously generates a Cartan connection by pull-back to $T \mathcal{P}$. From now on, we shall use Cartan and Tractor connections interchangeably.

\begin{rem}
It is not the case that any $G$ connection $\eta$ will correspond to a Cartan connection via pull-back to $\mathcal{P}$, as the isomorphism condition $T\mathcal{P}_v \to \mathfrak{g}$ could be violated. In fact, $\eta$ must have a maximal second fundamental form on the canonical sub-bundles in the splitting of the Tractor bundle. This form is sometimes known as the soldering form \cite{solder}. If so, then $\eta$ comes from a Cartan connection.
\end{rem}

\subsection{Conformal Geometry}

There are thus three standard ways of envisaging conformal geometry on a manifold $M$:
\begin{itemize}
\item[-] via a class of conformal metrics $[g]$ related by multiplication by a never-zero function (a zero order structure),
\item[-] via a class of torsion-free conformal connections $\nabla$ (a first order structure) or
\item[-] via a Cartan/Tractor Connection $\omega$/$\overrightarrow{\nabla}$ (a second order structure),
\end{itemize}

We will give more details of these three structures, and show their equivalence. The equivalence is easy to see between the first two structures -- $[g]$ defines a conformal frame bundle which is the principal bundle for the connections $\nabla$ -- but is non-trivial with the third structure.

Let $\mathfrak{co}(n) = \mathfrak{so}(n) \otimes \mathbb{R}$ be the conformal algebra. Then let $\mathcal{G}_0$ be the principal $\mathfrak{co}(n)$ frame bundle defined by $[g]$.

This allows us to define the bundles $\mathcal{E}[w]$, the weighted line-bundles coming from the centre of $\mathfrak{co}(n)$, i.e.
\be
\mathcal{E}[w] = \mathcal{G}_0 \times_{\rho} \mathbb{R} \ , \ \rho(c)(z) = -w \frac{\det(c)}{2n} z .
\ee
It is easy to see that $E[-n] = \wedge^{n} T^*$ and $E[w] = E[-n]^{-\frac{w}{n}}$, so any connection on the tangent bundle extends to a connection on these weighted line-bundles.

From now on, we'll use the notation $B[w]$ for $B \otimes \mathcal{E}[w]$. Then there is a map from any $g \in [g]$ to a section of $(\odot^2 T)[-2]$,
\be
g \to \mathbf{g} = \frac{n}{\det{g}} g.
\ee
This map does not depend on the choice of $g$. Conversely, given a non-vanishing section $\xi$ of $\mathcal{E}[1]$ -- a \emph{conformal scale} -- there is a corresponding metric in the conformal class
\be
g^{\xi} = \xi^{-2} \mathbf{g},
\ee
with a corresponding Levi-Civita connection $\nabla$. Thus the class $[g]$ and the \emph{conformal metric} $\mathbf{g}$ are equivalent, and we will use them interchangeably.

The second way of defining the conformal structure is to use the class of preferred connections:
\begin{defi}
Given a conformal manifold $\mathcal{G}_0 \to M$, a \emph{preferred connection} $\nabla$ is a torsion-free $\mathcal{G}_0$ connection.
\end{defi}
\begin{prop}
Given a conformal structure, a preferred connection is equivalent with a connection on $E[1]$ (or on any weighted line-bundle).
\end{prop}
\begin{proof}
Using the conformal metric to contract $T \otimes T \to E[2]$, we can use the same expression as for the Levi-Civita connection to generate a torsion-free conformal connection on the tangent bundle.
\end{proof}
In this view, those preferred connections that preserve a metric are exclusively those that preserve a conformal scale -- and hence have trivial curvature on $E[w]$.

\subsubsection{The Cartan connection}
In the classical, flat, case, Conformal Geometry is modelled on the sphere $S^n$. Taking the sphere as the collection of null-lines in $\mathbb{R}^{n+1,1}$, the group $G$ of conformal transformations is $SO(n+1,1)$. Then its Lie algebra has a 1-grading,
\begin{eqnarray*} \mathfrak{g} = \mathbb{R}^{n} \oplus \ \mathfrak{co}(n) \oplus \mathbb{R}^{n*},\end{eqnarray*}
where the conformal group $\mathfrak{co}(n)$ decomposes into the semisimple part $\mathfrak{so}(n)$ and the centre $\mathbb{R}$, which is responsible for the conformal weight in representations of $\mathfrak{co}(n)$.

Thus the data are $\mathfrak{g}=\mathfrak{so}(n+1,1)$, $\mathfrak{g}_0 = \mathfrak{co}(n)$ and $\mathfrak{p} = \mathfrak{co}(n) \rtimes \mathbb{R}^{n*}$, on an $n$-dimensional space. Note that we have a natural action of $\mathfrak{g}_0$ on $\mathfrak{g}$ and hence an associated bundle to the $\mathcal{G}_0$ structure bundle:
\begin{eqnarray*} \mathcal{G}_0 \times_{G_0} \mathfrak{g}. \end{eqnarray*}
Moreover, the action of $G_0$ splits $\mathfrak{g}$, giving a corresponding splitting:
\begin{eqnarray*} \mathcal{G}_0 \times_{G_0} \mathfrak{g} = T \oplus \mathfrak{co}(T) \oplus T^*. \end{eqnarray*}
This decomposition will be used extensively.

It is important to explicate the Lie bracket of this algebra \cite{TCPG}. In fact, $[T , T] = [T^*, T^*] = 0$, the Lie bracket on $\mathfrak{co}(T)$ is the natural commutator of endomorphisms, and $[c, t^*] = -c(t^*)$, $[c, t] = c(t)$, for $t, c, t^*$ sections of $T, \mathfrak{co}(T)$ and $T^*$ respectively. The bracket between $T$ and $T^*$ is more complicated, and in fact
\begin{eqnarray*} [t, s^*] = t \otimes s^* - (t \otimes s^*)^{\tau} + s^*(t) \delta, \end{eqnarray*}
with $\tau$ the transpose operator, and $\delta$ the identity element in (the centre of) $\mathfrak{co}(T)$.

In their papers \cite{TCPG} and \cite{TBIPG}, the authors demonstrate that the Cartan connection is equivalent to the standard conformal structure on a manifold $(M, [g])$. This is an alternative treatment.

\begin{theo} Let $s$ be a section of any bundle associated to $\mathcal{G}_0$, and let $X$ be any vector field. Then if $\nabla$ and $\widehat{\nabla}$ are two preferred connections, there exists a one-form $\Upsilon$ such that
\begin{eqnarray*} \widehat{\nabla}_X s = \nabla_X s + [\Upsilon, X]. s,\end{eqnarray*}
where $[,]$ is the Lie bracket for $\mathfrak{g}$ previously described.
\end{theo}
\begin{proof}
For $X \in \Gamma(T)$ and $\Upsilon \in \Omega^1(M)$, $[\Upsilon, X]$ is a section of $\mathfrak{co}(T)$, so this identity makes sense.

We know that $\widehat{\nabla}_X = \nabla_X + q(X)$, where $q$ is a one-form with values in $\mathfrak{co}(T)$.

However preferred connections are torsion-free, so $q(X).Y$ must be symmetric in $X$ and $Y$, implying that $q$ lies in the bundle $Q = \left( \odot^2 T^* \otimes T \right) \cap \left(T^* \otimes \mathfrak{co}(T) \right)$, the symmetrisation of $T^* \otimes \mathfrak{co}(T)$ around the first two elements.

The fact that $Q$ is of rank $n$ and spanned by elements of the form $[ \Upsilon, - ]$ can be seen by fairly simple Lie algebra manipulations (for more details, see \cite{methesis}; the important idea is that $Q \cap T^* \otimes \mathfrak{so}(T) =0$ by the uniqueness of the Levi-Civita connection, ensuring that the rank of $Q$ is $\leq n$).
\end{proof}

Note that if $\nabla$ were a metric connection, then $\widehat{\nabla}$ would be metric if and only if $\Upsilon$ were a closed form. In fact:
\begin{prop} Let $\nabla$ and $\widehat{\nabla} = \nabla + \Upsilon$ be two metric, preferred connections, with $\xi, \hat{\xi}$ the corresponding conformal scales. Defining the function $f$ as $f  = \xi \otimes \hat{\xi}^{-1} \in \Gamma(\mathcal{E}[0]) = C^{\infty}(M)$, we have
\begin{eqnarray*} \Upsilon = f^{-1} df = d (log f). \end{eqnarray*} \end{prop}
\begin{proof} By direct calculation, using the fact that $\nabla$ annihilates $g^{\xi}$ while $\widehat{\nabla}$ annihilates $g^{\hat{\xi}} = f^2 (g^{\xi})$.\end{proof}

A variety of tensors connected with these preferred connections will be needed in subsequent chapters. To define them, we will use Penrose's abstract index notation, where $Q^i$ is understood as a section of the tangent bundle, $Q_i$ a section of the co-tangent bundle, and symmetrisation and anti-symmetrisation of indexes to be denoted by $(ij)$ and $[ij]$ respectively. This notation will be used intermittently throughout the paper.

Then if $R_{ijkl}$ is the curvature tensor of $\nabla$, recall \cite{CCG}:
\begin{eqnarray} \label{W:def} R_{ijkl} = W_{ijkl} + 2 \mathbf{g}_{k[i}\mathsf{P}_{j]l} - 2 \mathbf{g}_{l[i}\mathsf{P}_{j]k} -2 \mathsf{P}_{[ij]}\mathbf{g}_{kl}\end{eqnarray}
with $W_{ijkl}$ the conformally invariant Weyl tensor, and the rho-tensor $\mathsf{P}$:
\begin{eqnarray} \label{confrho:formula} \mathsf{P}_{ij} = -\frac{1}{n-2} \ ( \frac{1}{n} \mathsf{Ric}_{ij} \ + \ \frac{n-1}{n} \mathsf{Ric}_{ji} \ - \ \frac{1}{2n-2} R \mathbf{g}_{ij}) \end{eqnarray}
a particularly important tensor for the rest of the paper. Here, $\mathsf{Ric}_{ij}$ is the Ricci curvature, and $R$ the scalar curvature $\mathsf{Ric}_{ij} \mathbf{g}^{ij}$ - a section of $\mathcal{E}[2]$.

This is in the general case for a conformal connection; in the metric case, the picture is the same, except that $\mathsf{P}$ follows the simpler symmetric formula
\begin{eqnarray*} \mathsf{P}_{ij} = - \frac{1}{n-2} \ ( \mathsf{Ric}_{ij} \ - \ \frac{1}{2n-2} R \mathbf{g}_{ij}). \end{eqnarray*}

The last relevant tensor for $\nabla$ is the Cotton-York tensor:
\begin{eqnarray} \label{CY:def} CY_{ijk} = 2 \nabla_{[i} \mathsf{P}_{j]k}.  \end{eqnarray}

It will be important to understand how the tensor $\mathsf{P}$ varies under a change of conformal structure, as this formula is the key to defining the Tractor bundle. Letting $\mathsf{P}$ be the rho-tensor for $\nabla$ and $\widehat{\mathsf{P}}$ be that of $\widehat{\nabla} = \nabla + \Upsilon$,
\begin{eqnarray} \label{p:change} \widehat{\mathsf{P}}(\xi) = \mathsf{P}(\xi) -\nabla_{\xi} \Upsilon + \frac{1}{2} [\Upsilon, [\Upsilon, \xi]], \end{eqnarray}
for $\xi$ any vector field.

\subsubsection{Equivalences}

Here we will demonstrate the equivalence of the Cartan connection with the conventional conformal structure. Though we will draw heavily on \cite{TBIPG} for this exposition, we will use a slightly unconventional approach, which has the advantage of constructing the vital `Tractor Bundle' directly.

\begin{rem} For a variety of reasons to do mainly with conventional notation and ease of calculations, we will be working with the Tractor bundle $\mathcal{T}$ in the rest of the paper. However, to get a better understanding of what this bundle actually is, we need to start by defining the dual bundle $\mathcal{T}^*$. \end{rem}

Consider the two-jet prolongation of $J^2(\mathcal{E}[1])$ of the weighted bundle $\mathcal{E}[1]$. By definition, we have the short exact sequences 
\begin{eqnarray*} 0 \longrightarrow \odot^2 T^*[1] \longrightarrow &J^2 (\mathcal{E}[1])& \longrightarrow J^1 (\mathcal{E}[1]) \longrightarrow 0, \\ 0 \longrightarrow T^*[1] \longrightarrow & J^1 (\mathcal{E}[1])& \longrightarrow \mathcal{E}[1] \longrightarrow 0. \end{eqnarray*}

The conformal structure $\mathbf{g}$ contracts $\odot^2 T^*$ to $\mathcal{E}[-2]$. Hence $\odot^2 T^*$ splits as $(\odot^2 T^*)_0 \oplus \mathcal{E}[-2]$, where the first space is the kernel of the contraction. Then we define the dual Tractor bundle $\mathcal{T}^*$ as the quotient:
\begin{eqnarray*} 0 \longrightarrow (\odot^2 T^*)_0 [1] \longrightarrow  & J^2(\mathcal{E}[1]) & \longrightarrow \mathcal{T}^* \longrightarrow 0. \end{eqnarray*}

It is actually possible to realise $\mathcal{T}^*$ as a sub-bundle of $J^2(\mathcal{E}[1])$ rather than a quotient bundle; we shall not be needing this result, though. Let $D$ be the second order operator $\Gamma(\mathcal{E}[1]) \to \Gamma(\mathcal{T}^*)$ given by composing the projection $J^2(\mathcal{E}[1]) \to \mathcal{T}^*$ with the two-jet operator $j^2$.

\begin{prop} \label{iso} Given a preferred connection $\nabla$, $s$ any section of $\mathcal{E}[1]$ and $b$ any point on the manifold, the map
\begin{eqnarray*} D s(b) \to (s(b), \nabla_i s (b), \frac{1}{n} \mathbf{g}^{ij} (- \nabla_i \nabla_j s(b) + \mathsf{P}_{ij} s(b))) \end{eqnarray*}
generates an isomorphism $\mathcal{T}^* \to \mathcal{E}[1] \oplus T^*[1] \oplus \mathcal{E}[-1]$.
\end{prop}

\begin{proof}
This formula clearly generates a bundle map $J^2(\mathcal{E}[1]) \to \mathcal{E}[1] \oplus T^*[1] \oplus \mathcal{E}[-1]$. All that remains is to prove that $(\odot^2 T^*)_0 [1]$ is the kernel of this map. Assume $D s (b) = 0$.

Then obviously $j^1(s) =0$ at $b$, implying that $\nabla_i \nabla_j s(b)$ is the (well-defined) section of $\odot^2 T^*$ that corresponds to the second derivative of $s$ at $b$. Thus $- \mathbf{g}^{ij} \nabla_i \nabla_j s(b) =0$, or equivalently $j^2(s)(b) = \nabla_i \nabla_j s(b) \in (\odot^2 T^*)_0$. \end{proof}

Notice that we have not used the tensor $\mathsf{P}$ yet; the next proposition shows what we need it for.

\begin{prop} Under a change of preferred connection $\nabla \to \widehat{\nabla} = \nabla + \Upsilon$, the isomorphism of Proposition \ref{iso} changes as
\begin{eqnarray*} \left( \begin{array}{c} x \\ \omega_i \\ z \end{array} \right) \to \left( \begin{array}{c} x \\ \omega_i - \Upsilon_i x \\ z + \mathbf{g}^{ij}\omega_i \Upsilon_j  - \frac{1}{2} \mathbf{g}^{ij} \Upsilon_i \Upsilon_j x \end{array} \right). \end{eqnarray*}
\end{prop}

\begin{proof} Direct computation from the transformation properties of $\nabla$ and $\mathsf{P}$. The first component obviously stays the same. As $\widehat{\nabla} = \nabla + \Upsilon$, we have the transformation law for the second component. For the third component:
\begin{eqnarray*} \widehat{\nabla}_i \widehat{\nabla}_j = (\nabla_i + \Upsilon_i) (\nabla_j + \Upsilon_j) = \nabla_i \nabla_j + \nabla_i \Upsilon_j + \Upsilon_i \nabla_j + \Upsilon^2_{ij}. \end{eqnarray*}
We need to remember that $\Upsilon_i$ acts on a section $x$ of $\mathcal{E}[1]$ as $\Upsilon_i . x = - \Upsilon_i x$.

Then recall the transformation law for $\mathsf{P}$ in equation $(\ref{p:change})$:
\begin{eqnarray*} \widehat{\mathsf{P}}_{ij} = \mathsf{P}_{ij} -\nabla_{j} \Upsilon_{i} + \frac{1}{2} \Upsilon^2_{ij}.  \end{eqnarray*}
Thus in total:
\begin{eqnarray*} (-\widehat{\nabla}_i \widehat{\nabla}_j + \widehat{\mathsf{P}}_{ij}).x &=& (-\nabla_{i} \nabla_{j} + \mathsf{P}_{ij}).x + \Upsilon_i  \nabla_j \ x - \frac{1}{2} \Upsilon^2_{ij} \ x \\ & & -(\nabla_{j} \Upsilon_{i} + \nabla_{i} \Upsilon_{j}) \ x \end{eqnarray*}
The last term in brackets is anti-symmetric, so disappears upon taking the contraction with the symmetric $\mathbf{g}^{ij}$. Hence, for $D(s) = (x, \omega_i, z)$ at $b$,
\begin{eqnarray*} \frac{1}{n} \mathbf{g}^{ij} (-\widehat{\nabla}_i \widehat{\nabla}_j + \widehat{\mathsf{P}}_{ij}) s &=& \frac{1}{n} \mathbf{g}^{ij} (-\nabla_{i} \nabla_{j} + \mathsf{P}_{ij} + \Upsilon_i \circ \nabla_j - \frac{1}{2} \Upsilon^2_{ij}) s  \\ &=& z + \mathbf{g}^{ij}\omega_i \Upsilon_j  - \frac{1}{2} \mathbf{g}^{ij} \Upsilon_i \Upsilon_j x. \end{eqnarray*}
\end{proof}

However, for reasons of convenience and notation, we will be working not with the bundle $\mathcal{T}^*$ but with its dual. Define the \emph{Tractor Bundle} as $\mathcal{T} = (\mathcal{T}^*)^*$.

The previous results carry through to the dual of $\mathcal{T}^*$; any preferred connection $\nabla$ defines a splitting $\mathcal{T} = \mathcal{E}[1] \oplus T[-1] \oplus \mathcal{E}[-1]$, and under a change of connection, this splitting changes via
\begin{eqnarray*} \left( \begin{array}{c} x \\ Y \\ z \end{array} \right) \to \left( \begin{array}{c} x \\ Y + \Upsilon^* x \\ z - \Upsilon(Y) - \frac{1}{2} \mathbf{g} (\Upsilon, \Upsilon) x \end{array} \right), \end{eqnarray*}
where $\Upsilon^* \in T[-2]$ is the dual to $\Upsilon$ using the conformal metric $\mathbf{g}$.

This particularly nice change of splitting formula implies the next vital lemma:
\begin{lemm} There is a natural metric $\langle , \rangle$, of type $(n+1,1)$, on $\mathcal{T}$. \end{lemm}
\begin{proof} Given a preferred connection $\nabla$ and two sections of $\mathcal{T}$, $(x, Y , z)$ and $(x', Y ' , z')$, we define the metric by
\begin{eqnarray*} \langle \left( \begin{array}{c} x \\ Y \\ z \end{array} \right) , \left( \begin{array}{c} x' \\ Y' \\ z' \end{array} \right) \rangle = xz' + x'z + \mathbf{g}(Y, Y'). \end{eqnarray*}
Direct calculation then shows this formula is invariant under a change of splitting. \end{proof}

Since $\mathcal{T}^*$ came about as a quotient bundle of a jet-bundle, we have invariant subspaces of $\mathcal{T}$:
\begin{eqnarray*} \mathcal{E}[-1] \ \hookrightarrow \ T \oplus \mathcal{E}[-1] \ \hookrightarrow \ \mathcal{T}, \end{eqnarray*}
and invariant projections:
\begin{eqnarray*} \pi^1 :& \mathcal{T} \to & \mathcal{E}[1] \oplus T, \\ \pi^2 :& \mathcal{T} \to & \mathcal{E}[1]. \end{eqnarray*}

Call $E$ the sub-bundle of $\mathcal{T}$ that is the inclusion of $\mathcal{E}[-1]$. Note that $E$ is null under $\langle,\rangle$.

\begin{lemm} There is a $\mathfrak{p}$-bundle $\mathcal{P}$ which is a principal bundle for $\mathcal{T}$. \end{lemm}
\begin{proof}
The metric $\langle , \rangle$ shows that the structure algebra of $\mathcal{T}$ reduces to $\mathfrak{g} = \mathfrak{so}(n+1,1)$. The invariant null sub-bundle $E$ further reduces the structure algebra to $\{ z \in \mathfrak{g} | z(E) = 0 \}$, i.e.~to $\mathfrak{p}$.

Then we define $\mathcal{P}$ to be the bundle of orthonormal frames of $\mathcal{T}$ preserving $E$.
\end{proof}

Let us review what we have so far. Starting from the conformal metric $\mathbf{g}$ and the class of preferred connections, we have constructed, via a tensor $\mathsf{P}$ dependent on the connections, a bundle $\mathcal{T}$. And this bundle generates a principal bundle $\mathcal{P}$, where it is natural to suppose the Cartan connection living. We now need to build this Cartan connection.

Define the Lie Algebra bundle $\mathcal{A} = \mathcal{P} \times_P \mathfrak{g}$. Then given a preferred connection $\nabla$ we have a splitting of $\mathcal{T} = \mathcal{E}[1] \oplus T[-1] \oplus \mathcal{E}[-1]$, and hence a corresponding splitting:
\begin{eqnarray*} \mathcal{A} = T \oplus \mathfrak{co}(T) \oplus T^*. \end{eqnarray*}

In order to finish the construction of this Cartan connection, we will start by building a $G = SO(n+1,1)$ connection on $\mathcal{A}$ and then prove that it is a Tractor connection.

\begin{defi} Given a preferred connection $\nabla$, we have a splitting of $\mathcal{A} = T \oplus \mathfrak{co}(T) \oplus T^*$. Each of these bundles is a $G_0$ bundle, so $\nabla$ ascends to a connection on $\mathcal{A}$. Then we define the Tractor connection $\overrightarrow{\nabla}$ as
\begin{eqnarray*} \overrightarrow{\nabla}_X = \nabla_X + \mathrm{ad}(X) + \mathrm{ad} \mathsf{P}(X),
\end{eqnarray*}
with the vector $X$ and the one-form $\mathsf{P}(X)$ seen as sections of the Lie algebra bundle $\mathcal{A}$.
\end{defi}

Then since $\nabla$ is a $\mathfrak{g}_0$ connection, $\overrightarrow{\nabla}$ is a $\mathfrak{g}$ connection. Of course, this definition makes no sense without:
\begin{prop} This definition is independent of the choice of $\nabla$. \end{prop}
\begin{proof} The formula for the change of splitting of $\mathcal{A}$ (deduced directly from that of $\mathcal{T}$) is:

\begin{eqnarray*} \left( \begin{array}{c} X \\ \Psi \\ \omega \end{array} \right) \to \left( \begin{array}{c} X \\ \Psi + [\Upsilon, X] \\ \omega + [\Upsilon, \Psi] + \frac{1}{2} [\Upsilon, [\Upsilon, X]] \end{array} \right), \end{eqnarray*}
Then a direct calculation proves the result.
\end{proof}

Thus for any bundle $B$ associated to $\mathcal{G}$, we have an invariant connection form:
\begin{eqnarray*} \overrightarrow{\nabla} = \nabla + \rho(X) + \rho \mathsf{P}(X). \end{eqnarray*}

In the case of the Tractor bundle $\mathcal{T}$, the detailed expression is:
\begin{eqnarray} \label{basic:eqn} \overrightarrow{\nabla}_{X} \left( \begin{array}{c} x \\ Y \\ z \end{array} \right) = \left( \begin{array}{c} \nabla_{X} x - \mathbf{g}(X, Y) \\ \nabla_{X} Y + z X - x \mathsf{P}(X) \\ \nabla_{X} z + \mathsf{P}(X, Y) \end{array} \right), \end{eqnarray}

Now we get to the result that ties all the structures together:
\begin{theo} The connection \overrightarrow{\nabla} is a normal Tractor connection. \end{theo}
\begin{proof}
Let $i$ be the inclusion $i: \mathcal{P} \hookrightarrow \mathcal{G}$, $\pi$ projection $\pi: \mathcal{P} \to M$, and let $\omega \in \Omega^1(\mathcal{G}, \mathfrak{g})$ be the one-form associated with the connection $\overrightarrow{\nabla}$.

We need to prove that $\mu = i^*(\omega)$ is an isomorphism $T \mathcal{P}_{u} \to \mathfrak{g}$ for all points $u \in \mathcal{P}$; then $\mu$ will be the Cartan connection generating the Tractor connection $\omega$.

So now assume that $\mu$ is \emph{not} an isomorphism at some point $u$, so there exists a vector $\xi \in T \mathcal{P}_u$ such that $\mu(\xi) = 0$. As vertical vectors in $\mathcal{P}$ are mapped isomorphically onto $\mathfrak{p}$, $X = \pi_* (\xi)$ is a non zero vector in $T_{\pi(u)}$.

Then define a local section $j$ of $M$ in $\mathcal{P}$, such that $j_* X = \xi$. This also gives us a section $i \circ j$ of $\mathcal{G}$. Then in the frame bundle determined by this section, the connection $\overrightarrow{\nabla}$ is of the form
\begin{eqnarray*} d + j^* i^* \omega = d + j^* \mu \end{eqnarray*}

This shows that the bundle $\mathcal{P}$ is infinitesimally conserved at $\pi(u)$ in the $X$ direction, or, switching to the associated bundle $E = \mathcal{P} \times_P e$, that there is a section $s$ of $E$, non-zero at $\pi(u)$, such that $\overrightarrow{\nabla}_X s =0 $ at $\pi(u)$.

However, the connection on the tractor bundle is given by equation (\ref{basic:eqn}):
\begin{eqnarray*} \overrightarrow{\nabla}_{X} \left( \begin{array}{c} 0 \\ 0 \\ z \end{array} \right) = \left( \begin{array}{c} 0 \\ z X \\ \nabla_{X} z \end{array} \right), \end{eqnarray*}
which is a contradiction as $s$ (hence $z$) is non-zero at $\pi(u)$. So $\overrightarrow{\nabla}$ is indeed a Tractor connection.

And finally, to complete the circle:
\begin{lemm} \label{norm:lemm} The Cartan connection generated by $\overrightarrow{\nabla}$ is normal. \end{lemm}
\begin{proof}
By \cite{TCPG} and \cite{kos}, this result is equivalent with the curvature of $\overrightarrow{\nabla}$ lying in the Lie algebra bundle of $\mathcal{G} \times_G \mathfrak{p}$. Alternately, the curvature must preserves the canonical bundle $E$.

In abstract index notation, the expression for $\overrightarrow{\nabla}_i \overrightarrow{\nabla}_j$ is:
\begin{eqnarray*} \overrightarrow{\nabla}_i \overrightarrow{\nabla}_j &=& \nabla_i \nabla_j + (\delta^l_j \nabla_i + \delta^l_i \nabla_j) +  (\mathsf{P}_{jk} \nabla_i + \mathsf{P}_{ik} \nabla_j) + \nabla_i (\mathsf{P}_{jk}) \\ & & + \delta_i^l \circ \delta_j^m + \delta_i^l \circ \mathsf{P}_{jn} + \mathsf{P}_{ik} \circ \delta_j^m + \mathsf{P}_{ik} \circ \mathsf{P}_{jn}. \end{eqnarray*}
Here we have used the connection $\nabla$ on $T$ to define the second covariant derivative; however, we could have used any other connection, as we are about to anti-symmetrise $i$ and $j$. Upon doing this, the terms in brackets vanish. Moreover, $\rho(X) \circ \rho(Y) = \rho(Y) \circ \rho(X)$ and similarly for one-forms, meaning that:
\begin{eqnarray*} R^{\overrightarrow{\nabla}}_{ij} &=& R^{\nabla}_{ij} + \nabla_{[i} (\mathsf{P}_{j]k}) + \delta_{[i}^l \circ \mathsf{P}_{j]n} + \mathsf{P}_{k[i} \circ \delta_{j]}^m. \end{eqnarray*}

Looking back at equations ($\ref{W:def}$) and ($\ref{CY:def}$), we see that this expression is the sum of the Weyl tensor and the Cotton York tensor. Or, expressed in more conventional notation, in the splitting of $\mathcal{A}$ determined by $\nabla$:
\begin{eqnarray*} \label{curv:def} R^{\overrightarrow{\nabla}}_{X,Y} = \left( \begin{array}{c} 0 \\ W(X,Y) \\ CY(X,Y) \end{array} \right). \end{eqnarray*}
Since only $T \subset \mathcal{A}$ has a non-trivial action on the canonical bundle $E \subset \mathcal{T}$, this curvature expression must preserve $E$.
\end{proof}
\end{proof}
To construct the conformal structure from the Cartan connection is much simpler; indeed, the metric $\langle, \rangle$ descendes to the conformal metric $\mathbf{g}$ on $T = E^{\perp}/E$.

\subsubsection{One and two dimensions} \label{two:dim}
Though any two-manifold is conformally flat, with an infinite-dimensional local conformal tranformation group, paper \cite{Mob} and other unpublished papers by the same author extend the concept of conformal Cartan connections to one and two dimensions, by constructing M\"obius structures. As in higher dimensions, a choice of Weyl structure determines a splitting of the associated Tractor bundle. There is an ambiguity, however, in the trace-free symmetric part of the $\mathsf{P}$-tensor; this may be chosen freely.

\begin{defi}
For our purposes, we shall take
\begin{eqnarray*}
\mathsf{P}_{hj} &=& -\frac{1}{2} \mathsf{Ric}_{jh}.
\end{eqnarray*}
This is not a conformally invariant definition. However, we shall be using it in a specific metric (Einstein, with constant scalar curvature), where it makes sense and allows one to extend the reach of the decomposition theorem down to lower dimensions.
\end{defi}

In one dimension, we may easily require
\begin{eqnarray*}
\mathsf{P}_{hj} = 0,
\end{eqnarray*}
which is conformally invariant. This also fits our definitions.

\section{Conformally Einstein Manifolds}

\subsection{Important Note}

In most of the proofs in the remainder of this paper, it will be assumed that for a certain holonomy preserved sub-bundle $U \subset \mathcal{T}$ used in the proof, one has $\pi^2(U) \neq 0$. This will not be the case everywhere, of course; however:

\begin{prop} \label{HOL:LOC} Let $U \subset \mathcal{T}$ be a preserved subbundle under $\overrightarrow{\nabla}$. Then $\pi^2(U) \neq 0$ on $\Sigma$, an open, dense subset of $M$. \end{prop}
\begin{proof}

$\Sigma$ is open because of the $\pi^2(U) \neq 0$ condition.

Let $b \in M \backslash \Sigma$, and $u(b) = (0, Y(b), z(b))$ be a non-zero element of $U_b$. Then extend $u(b)$ locally to a section $u = (x, Y, z)$ of $U$ by parallel transport along `rays' from $b$. This implies that $\overrightarrow{\nabla} u = 0$ at $b$. Then picking any nowhere-zero section $\tau$ of $\mathcal{E}[-1]$, we can define the function $f: M \to \mathbb{R}^n$ by
\begin{eqnarray*} f(c)= \pi^2(u) \tau. \end{eqnarray*}

The derivative of $f$ is $(\nabla x) \tau + x (\nabla \tau)$. At $b$, this is just $(\nabla x ) \tau$, and, since $\overrightarrow{\nabla} u = 0$ at $b$:
\begin{eqnarray*} X.f(b)= (\nabla_X x) \tau = - \mathbf{g}(X,Y) \tau. \end{eqnarray*}

If $Y \neq 0$, this is non-zero for some $X$, so $f$ is non-zero arbitrarily close to $b$. If $Y=0$, then the first derivative is zero, and the second derivative is thus:
\begin{eqnarray*} Z.(X.f)(b) &=& (\nabla_Z \nabla_X x) \tau \\ &=& - \mathbf{g}(X, \nabla_Z Y) \tau \\ &=& -\mathbf{g}(X, -z Z) \tau. \end{eqnarray*}
with $z(b) \neq 0$ as $u(b) \neq 0$. Then the second derivative is non-zero for $X = Z \neq 0$, for instance, forcing $f$ to be non-zero arbitrarily close to $b$.

This implies that $\pi^2(u) \neq 0$ arbitrarily close to $b$, proving the result. \end{proof}
In fact, if the first derivative vanishes, $b$ must be an \emph{isolated} point.

The classic examples of this are the various conformally Einstein metrics on the sphere $S^n$. The sphere is conformally flat, so there are many holonomy preserved sections of its Tractor bundle.

A preserved section $u$ of negative norm corresponds to the Spherical metric $g = \pi^2(u)^{-2} \mathbf{g}$ on the whole space. In this case, $\pi^2(u)$ is never zero.

A preserved section $u$ of zero norm corresponds to the Euclidean metric $g = \pi^2(u)^{-2} \mathbf{g}$ on $\mathbb{R}^n \cong S^n \backslash \{ \infty \}$. In this case, $\pi^2(u)(b) \neq 0$ for $b \neq \infty$.

A preserved section $u$ of positive norm corresponds to the Hyperbolic metric $g = \pi^2(u)^{-2} \mathbf{g}$ on two half spheres of $S^n$. In this case $\pi^2(u)$ is zero only on the boundary $S^{n-1}$ cutting $S^n$ into two.

\subsection{Einstein Spaces}

Though it is well known in general that any conformally Einstein space corresponds to a parallel section of the tractor bundle $\mathcal{T}$, what follows is a direct proof of this fact using the Tractor connection approach.

\begin{rem} This is a first instance of a holonomy reduction of $\overrightarrow{\nabla}$. \end{rem}

\begin{theo} For $n>2$, if $(M, \mathbf{g})$ has an Einstein metric $g$ in its conformal class then there exists a parallel section $s$ of its tractor bundle $\mathcal{T}$. \end{theo}

\begin{proof}

Let $g$ be the Einstein metric, $\mathsf{Ric}^g = \lambda g$. Then the $\mathsf{P}$-tensor is
\begin{eqnarray*} \mathsf{P} = - \frac{\lambda}{2n-2} \ \mathbf{g}. \end{eqnarray*}

Hence
\begin{eqnarray*} \overrightarrow{\nabla} \left( \begin{array}{c} 1 \\ 0 \\ - \frac{\lambda}{2n-2} \end{array} \right) = \left( \begin{array}{c} 0 \\ \frac{\lambda}{2n-2} - \frac{\lambda}{2n-2} \\ 0 \end{array} \right) = 0, \end{eqnarray*}
where $1$ is the section of $\mathcal{E}[1]$ corresponding to $g$.
\end{proof}

To prove the converse of this theorem, we need the following lemma:
\begin{lemm} If a conformal connection $\nabla$ has a symmetric Ricci tensor, then $\nabla$ is actually a metric connection. \end{lemm}

\begin{proof}
Let $R^i_{\phantom{i}jkl}$ be the curvature of $\nabla$. Then $R^i_{\phantom{i}jkl}$ acts on the determinant bundle $\mathcal{E}[n]$ via its trace $R^i_{\phantom{i}ikl}$. However, by the first Bianci identity,
\begin{eqnarray*}
R^i_{\phantom{i}ikl} &=& - R^i_{\phantom{i}kli} - R^i_{\phantom{i}lik} \\ &=& \mathsf{Ric}_{kl} - \mathsf{Ric}_{lk},
\end{eqnarray*}
the anti-symmetric part of the Ricci tensor. So if $\nabla$ has a symmetric Ricci tensor, its curvature must vanish on $\mathcal{E}[n]$, so locally $\nabla$ must preserve a section $\eta$ of the determinant bundle. Then $\eta^{1/n}$ is a preserved conformal scale and
\begin{eqnarray*} g = \eta^{-2/n} \mathbf{g} \end{eqnarray*}
a metric preserved by $\nabla$. \end{proof}

\begin{theo} \label{ric:flat} For $n>2$, if a line bundle $L$ of $\mathcal{T}$ is holonomy preserved, then a section $s$ of $L$ is preserved, and $(M, [g])$ has an Einstein metric $g = (\pi^2(s))^{-2} \mathbf{g}$ in its conformal class, wherever $\pi^2 (s) \neq 0$. \end{theo}

\begin{proof}

The line bundle $L$ defines a connection on $\mathcal{E}[1]$, and hence a torsion free connection on $T$, in the following way. Let $e$ be any nowhere vanishing section of $\mathcal{E}[1]$, and let $l$ be the section of $L$ such that $\pi^2(l) = e$. Then define $\nabla e = \pi^2(\overrightarrow{\nabla}l)$; it is easy to see that this is indeed a connection.

Using $\nabla$, we split $\mathcal{T} = \mathcal{E}[1] \oplus T[-1] \oplus \mathcal{E}[-1]$ in the usual way. Then equation (\ref{basic:eqn}) implies that $\pi^2(\overrightarrow{\nabla} (x, Y, z)) = \nabla x - \mathbf{g}(Y, - )$. Since by definition of $\pi^2 (\overrightarrow{\nabla} l) = \nabla \pi^2(l)$, we must have $Y = 0$ for $l$ (and hence for any section of $L$).

If $L$ is not null, then a section $s = (x, 0, x\mu)$ of constant norm, is preserved. This generates a metric $g = (\pi^2(s))^{-2} \mathbf{g}$. But $\overrightarrow{\nabla}_X (x, 0, x \mu)= (\nabla_X x, x \mu(X) -x\mathsf{P}(X), \mu \nabla_X x) = 0$. This implies that $\mathsf{P} = \mu g$, so $\mathsf{Ric}^{g} = \lambda g$ for $\lambda = (2-2n) \mu$.

On the other hand, if $L$ is null, $z=0$, and $\overrightarrow{\nabla}_X (x, 0, 0)= (\nabla_X x, -x\mathsf{P}(X), 0)$. Thus $\mathsf{P} = 0$ and hence $\nabla$ has a symmetric Ricci tensor, implying that it is actually a metric connection for some metric $g$ -- which moreover is Ricci-flat. Set $x \in \Gamma(\mathcal{E}[1])$ to be the conformal scale corresponding to g. Then the section $s = (x,0,0)$, is parallel.
\end{proof}

\begin{prop}
In the two dimensional case, one merely has the one-way implication that an Einstein metric of constant scalar curvature gives a preserved section of $\mathcal{T}$.
\end{prop}
\begin{rem} Note that the sign of $\langle s,s \rangle$ is the opposite of the sign of the Einstein constant $\lambda$. \end{rem}

\section{Decomposition Theorem}

This section presents the decomposition theorem for Tractor connections, similar to the De Rham decomposition for Riemannian connections.

\begin{rem} Related terminology may be found in $\cite{Nets}$. \end{rem}

\subsection{Preparatory Results}
\begin{defi}
Given a metric $g$ on $M$ with Levi-Civita connection $\nabla$, a subbundle $U \subset T$ is \emph{umbilical} for the connection $\nabla$, if there exists a vector field $H$ such that for $X$ and $Y$ sections of $U$,
\begin{eqnarray*}
\nabla_X Y = \widetilde{\nabla}_X Y + g(X,Y) H,
\end{eqnarray*}
for $\widetilde{\nabla}$ some connection on $U$, and $H$ a vector field.
\end{defi}
\begin{rem}
Note that an umbilical subbundle is automatically integrable, as
\begin{eqnarray*}
[X,Y] = \nabla_X Y -\nabla_Y X &=& \widetilde{\nabla}_X Y - \widetilde{\nabla}_Y X + \big( g(X,Y) - g(Y,X) \big) H \\
&=& \widetilde{\nabla}_X Y - \widetilde{\nabla}_Y X,
\end{eqnarray*}
a section of $U$.
\end{rem}
\begin{lemm}
$U$ being umbilical is equivalent to
\begin{eqnarray} \label{umb:equ}
\nabla_X Y & \in & \Gamma(U),
\end{eqnarray}
whenever $X$ and $Y$ are orthogonal sections of $U$.
\end{lemm}
\begin{proof}
If $U$ is umbilical, then Equation (\ref{umb:equ}) is true by definition
\begin{eqnarray*}
\nabla_X Y &=& \widetilde{\nabla}_X Y + g(X,Y) H \\
&=& \widetilde{\nabla}_X Y \in \Gamma(U).
\end{eqnarray*}
So we now assume Equation (\ref{umb:equ}) and aim to prove umbilicity. One may easily see, by choosing an orthogonal frame for $U$, that $U$ must be integrable.

Define a connection $\widetilde{\nabla}$ on $U$, by orthogonal projection. Then the map $\Phi = \nabla - \widetilde{\nabla}$ is bilinear, $U^* \otimes U^* \to U^{\perp}$, and symmetric since $\nabla_X Y - \nabla_Y X = [X,Y]$ is a section of $U$. By assumption, $\Phi(X,Y) = 0$ whenever $g(X,Y) = 0$.

Now let $(X_j)$ be a frame of $U$, chosen so that the $g(X_j, X_k)$ are nowhere zero (one can do this, for instance, by choosing a standard orthonormal frame $(X_j)$ and mapping $X_j \to X_j + \frac{1}{2n} \sum_{l=1}^j X_l$). Pick $H$ in $U^{\perp}$ such that $\Phi(X_1, X_1) = g(X_1, X_1)H$. Then since $X_1$ is orthogonal to $\tau_{1,1,j} = g(X_1, X_1) X_j - g(X_j, X_1)X_1$, one has $\Phi(X_1, \tau) = 0$ and hence
\begin{eqnarray*}
\Phi(X_1, X_j) &=& \frac{1}{g(X_1, X_1)} g(X_j, X_1) \big( g(X_1, X_1)H \big) \\
&=& g(X_j, X_1)H.
\end{eqnarray*}
The same argument with the orthogonal sections $\tau_{j,1,k}$ and $X_j$ demonstrates
\begin{eqnarray*}
\Phi(X_j, X_k) &=& g(X_j, X_k)H.
\end{eqnarray*}
This extends trivially to the whole of $U$. Thus $\nabla_X Y = \widetilde{\nabla}_X Y + \Phi(X,Y) = \widetilde{\nabla}_X Y + g(X,Y)H$.
\end{proof}

Note that being umbilical is a conformally invariant condition, as changing $\nabla$ by $\Upsilon$ changes $H$ to $H - \Upsilon^*$. Thus choosing $\Upsilon^* = H$, we can make $U$ into a totally geodesic foliation. In other words, there are preferred connections for which $U$ is totally geodesic.

\subsection{Preserved subbundles}

Let $K$ be a subbundle of $\mathcal{T}$ of rank $k$, $2 \geq k \leq n$, preserved by $\overrightarrow{\nabla}$. Then $K$ defines a sub-bundle $U$ of $T$ as follows. We assume, from Lemma \ref{HOL:LOC}, that $K$ and $K^{\perp}$ are locally transverse to $E$. Recall that $E \cong E[-1] \subset \mathcal{T}$ is the canonical line bundle, and that $E^{\perp} \cong TM[-1] \oplus L^{-1}$ is of rank $n+1$ in $\mathcal{T}$.

Hence, $K \cap E^{\perp}$ is a bundle of rank $k-1$, and $\pi^1$ is injective on $(K \cap E^{\perp})$ (since $K \cap E = 0$, so $\pi^1$ is injective on $K$). Moreover $\pi^1(E^{\perp}) = T[-1]$, so
\begin{eqnarray*}
U = \pi^1(K \cap E^{\perp}) \subset T[-1]
\end{eqnarray*}
is a well defined, rank $k-1$ bundle. Use any conformal scale to get an isomorphism $T \cong T[-1]$. Since changing the section simply results in scaling any element of $T[-1]$, we may see $U$ as a well-defined subbundle of $T$.

\begin{prop}
$U$ is an integrable, umbilical foliation of $T$.
\end{prop}
\begin{proof}
Let $X$ and $Y$ be orthogonal sections of $U$. Fix any metric in the conformal class. Then
\begin{eqnarray*}
\left( \begin{array}{c} 0 \\ Y \\ z \end{array} \right),
\end{eqnarray*}
is a section of $K \cap E^{\perp}$, for some $z$. Consequently
\begin{eqnarray*}
\overrightarrow{\nabla}_X \left( \begin{array}{c} 0 \\ Y \\ z \end{array} \right) = \left( \begin{array}{c} 0 \\ \nabla_X Y + zX \\ z' \end{array} \right),
\end{eqnarray*}
for some $z'$. Since $K$ is preserved by $\overrightarrow{\nabla}$, this is a section of $K$; it is clearly a section of $E^{\perp}$. As a consequence, we know that
\begin{eqnarray*}
\nabla_X Y + zX & \in & \Gamma(U).
\end{eqnarray*}
Thus $\nabla_X Y$ is also a section of $\Gamma(U)$, making $U$ umbilical, and hence integrable.
\end{proof}
We shall see later that $U$ is Einstein (i.e.~all leaves $N$ of $U$ are conformally Einstein under the restricted conformal structure).

\begin{prop}
There is a Tractor bundle $\mathcal{T}_U$ on the leaves $N$ of the foliation defined by $U$, and a well-defined inclusion $\mathcal{T}_U \subset \mathcal{T}$.
\end{prop}
\begin{proof}
If $\nabla$ is a $U$-preferred connection -- one that makes $U$, and its foliation, totally geodesic -- in the splitting of $\mathcal{T}$ that it defines,
\begin{eqnarray*}
\mathcal{T} = L^1 \oplus T[-1] \oplus L^{-1}.
\end{eqnarray*}
Define $\mathcal{T}_U$ as the subbundle
\begin{eqnarray*}
\mathcal{T}_U = L^1 \oplus U[-1] \oplus L^{-1}.
\end{eqnarray*}
To check this is well defined, we change $\nabla$ to $\nabla'$, another $U$-preferred connection. This is equivalent to changing $\nabla$ by an $\Upsilon \in \Gamma(g(U) \subset T^*)$ for any metric $g$ in the conformal class. Then the splitting changes as:
\begin{eqnarray*}
\left( \begin{array}{c} x \\ Y \\ z \end{array} \right) \to \left( \begin{array}{c} x \\ Y + \Upsilon^* x \\ z - \Upsilon(Y) - \frac{1}{2} \mathbf{g} (\Upsilon, \Upsilon) x \end{array} \right),
\end{eqnarray*}
which, since $\Upsilon^*$ is a section of $U$, does not change the definition of $\mathcal{T}_U$ nor its inclusion into $\mathcal{T}$.
\end{proof}
We are now ready to prove the main theorem.
\begin{theo} \label{split:theorem}
Assume there is a bundle $K$ of rank $k$ preserved by $\overrightarrow{\nabla}$, and the foliation $U$ that it generates splits $T$. Let $l = k-1$ be the rank of $U$. Then there exists a metric $g$ in the conformal class of $M$ such that the manifold $(M,g)$ splits locally as the direct product
\begin{eqnarray*}
(M,g) = (N_1, h_1) \times (N_2, h_2)
\end{eqnarray*}
where $h_1$ and $h_2$ are Einstein metrics with Einstein coefficients $\lambda_1$, $\lambda_2$, possibly zero, related by
\begin{eqnarray*}
(n-l-1) \lambda_1 = (1-l) \lambda_2.
\end{eqnarray*}
The converse is also true. And in this situation
the holonomy $\overrightarrow{\mathfrak{hol}}$ of $\overrightarrow{\nabla}$ is the direct sum of Lie algebras
\begin{eqnarray*}
\overrightarrow{\mathfrak{hol}} = \overrightarrow{\mathfrak{hol}}_{N_1} \oplus \overrightarrow{\mathfrak{hol}}_{N_2}
\end{eqnarray*}
where $\overrightarrow{\mathfrak{hol}}_{N_1}$ is the holonomy of $\overrightarrow{\nabla}_{N_1}$ and $\overrightarrow{\mathfrak{hol}}_{N_2}$ that of $\overrightarrow{\nabla}_{N_1}$.
\end{theo}

Note that the subbundle of $T$ generated by $K^{\perp}$ is just $U^{\perp}$. There are really two situations here: the case when $K \cap K^{\perp}$ is of rank one, and that where it is of rank zero.

\subsubsection{$K$ degenerate}
If $K \cap K^{\perp} = \mathcal{L}$, a line bundle, necessarily null, then by Theorem \ref{ric:flat} there must be a preserved section $v$ of $\mathcal{L}$ and hence a Ricci-flat metric $g$ on $M$, with Levi-Civita connection $\nabla$.

We have the bundles $U$ and $U^{\perp}$ as before, both integrable and umbilical. We will now show that $g$ is locally a product metric of the leaves genereated by $U$ and $U^{\perp}$. First, we shall demonstrate that these leaves are totally geodesic under $g$.
\begin{lemm}
Let $X$ be a section of $U$. Then for any $A \in \Gamma(T)$, $\nabla_A X$ is a section of $U$.
\end{lemm}
\begin{proof}
In the splitting defined by $g$, one section of $K$ is the Einstein vector
\begin{eqnarray*}
v = \left( \begin{array}{c} 1 \\ 0 \\ 0 \end{array} \right).
\end{eqnarray*}
Since $v$ is also a section of $K^{\perp}$, $K$ must lie in $v^{\perp}$. In other words, $K$ is of the form
\begin{eqnarray*}
\left( \begin{array}{c} \mathbb{R} \\ U \\ 0 \end{array} \right).
\end{eqnarray*}
Now consider
\begin{eqnarray*}
\overrightarrow{\nabla}_A \left( \begin{array}{c} 0 \\ X \\ 0 \end{array} \right) = \left( \begin{array}{c} -\mathbf{g}(A,X) \\ \nabla_A X \\ 0 \end{array} \right).
\end{eqnarray*}
Since $\overrightarrow{\nabla}$ preserves $K$, $\nabla_A X$ must be a section of $U$.
\end{proof}
This shows that $U$ (and $U^{\perp}$) are totally geodesic foliations. Moreover, they are preserved by $\nabla$ in every direction.
\begin{rem}
As a consequence of that, if $X$ and $B$ are commuting sections of $U$ and $U^{\perp}$ respectively,
\begin{eqnarray*}
\nabla_X B = \nabla_B X = 0.
\end{eqnarray*}
\end{rem}

Let $h_1 = g|_U$, $Y$ and $X$ be sections of $U$, $A$ any section of $T$. Then
\begin{eqnarray*}
(\nabla_A h_1)(X,Y) &=& A. h_1(X,Y) - h_1 (\nabla_A X, Y) - h_1(X, \nabla_A Y) \\
&=& A. g(X,Y) - g (\nabla_A X, Y) - g(X, \nabla_A Y) \\
&=& (\nabla_A g)(X,Y) \\
&=& 0,
\end{eqnarray*}
as $\nabla_A X$ and $\nabla_A Y$ are sections of $U$, and $h_1 = g$ on sections of $U$. Consequently we have demonstrated, for $h_1$ and for $h_2 = g|_{U^{\perp}}$:
\begin{lemm}
$\nabla h_1$ and $\nabla h_2$ are both zero.
\end{lemm}
Now pick sections $X$ and $Y$ of $U$ commuting with a section $B$ of $U^{\perp}$. By the previous lemma
\begin{eqnarray*}
B. h_1(X,Y) = 0,
\end{eqnarray*}
so the Lie derivative of $h_1$ in the direction of $B$ is
\begin{eqnarray*}
\big( \mathcal{L}_B h_1 \big) (X,Y) = B. h_1(X,Y) - h_1([B,X],Y) - h_1(X,[B,Y]) = 0.
\end{eqnarray*}
We may choose local coordinates that respect the foliations $U$ and $U^{\perp}$ to get frames $(X^j)$ of $U$ and $(B^k)$ of $U^{\perp}$, commuting with one-another. Consequently, if $N_1$ is a leaf of $U$ and $N_2$ a leaf of $U^{\perp}$, $h_1$ is preserved by translation along $N_2$ and vice-versa. This demonstrates that
\begin{prop} \label{times:manifold}
Locally, $(M,g) = (N_1,h_1) \times (N_2, h_2)$.
\end{prop}
This implies that $\nabla|_U$ is the Levi-Civita connection of $h_1$, and $\nabla|_{U^{\perp}}$ that of $h_2$. To finish this exploration, we require:
\begin{lemm}[Restricted Ricci curvature] \label{ricci:rest}
Given a foliation $U$ preserved by $\nabla$, the Ricci tensor of $\nabla|_U$ is the Ricci tensor of $\nabla$, restricted to $U$.
\end{lemm}
\begin{proof}
Notice that this condition makes $U$ integrable and totally geodesic. Let $(X_j),(B_j)$ be a coordinate frame for $T$, with $X_j \in \Gamma(U)$ and the $(B_j)$ complementary. Then
\begin{eqnarray*}
\mathsf{Ric}(X_j,X_k) = \left( \sum_l X_l^* \lrcorner R_{X_l, X_j} X_k \right) + \left( \sum_l B_l^* \lrcorner R_{B_l, X_j} X_k \right).
\end{eqnarray*}
But the second term on the right is zero, as $R_{-,-} X_j$ must be a section of $U$, and the first term is just the Ricci curvature of $\nabla_U$.
\end{proof}
Consequently, one can see that $\nabla$ is Ricci-flat on $U$ and on $U^{\perp}$ (hence on $N_1$ and $N_2$).

In this case the relation
\begin{eqnarray*}
(n-l-1) \lambda_1 = (1-l) \lambda_2.
\end{eqnarray*}
is trivially satisfied, as both $\lambda_j$ are zero. The converse to this construction is trivial: a direct product of Ricci-flat spaces is Ricci-flat. Then $K$ may be reconstructed as
\begin{eqnarray*}
K = \left( \begin{array}{c} \mathbb{R} \\ TN_1 \\ 0 \end{array} \right)
\end{eqnarray*}
in the global Ricci-flat metric's splitting. Since $TN_1$ must be totally geodesic, $\overrightarrow{\nabla}$ preserves $K$ and
\begin{eqnarray*}
K^{\perp} = \left( \begin{array}{c} \mathbb{R} \\ TN_2 \\ 0 \end{array} \right).
\end{eqnarray*}
Now notice that since all $\mathsf{P}$ are zero, $\overrightarrow{\nabla}$ acts on $\mathcal{T}_{N_1}$ along $N_1$ exactly as the Tractor connection $\overrightarrow{\nabla}_{N_1}$ does. Moreover, $\overrightarrow{\nabla}$ acts trivially on $\mathcal{T}_{N_1}$ along $N_2$. Since the opposite result holds for $\mathcal{T}_{N_2}$, and since these two tractor bundles span all of $\mathcal{T}$, one has
\begin{eqnarray*}
\overrightarrow{\mathfrak{hol}} = \overrightarrow{\mathfrak{hol}}_{N_1} \oplus \overrightarrow{\mathfrak{hol}}_{N_2}.
\end{eqnarray*}

\subsubsection{$K$ non-degenerate}
We seek to imitate the proofs of the previous section in the case where $K \cap K^{\perp} = 0$. First of all, we seek to find an imitation of the Ricci-flat metric $g$. We shall use a preferred connection rather than a metric -- though it will turn out to be a metric connection in the end.

Starting off, pick $\nabla'$ such that $U$ is totally geodesic. In the rest of these proofs, $X$ and $Y$ will be sections of $U$, $B$ and $C$ sections of $U^{\perp}$.

Since $U^{\perp}$ is umbilical,
\begin{eqnarray*}
\nabla'_B C = \widetilde{\nabla}'_B C + H \tilde{g} (B,C),
\end{eqnarray*}
for some $H \in \Gamma(U)$ and any metric $\tilde{g}$ in the conformal class. Then replace $\nabla'$ with $\nabla$, by adding the one-form $\Upsilon = \tilde{g}(H)$. This connection makes $U^{\perp}$ totally geodesic, but since
\begin{eqnarray*}
\nabla_X Y = \nabla'_X Y + \Upsilon(X) Y + \Upsilon(Y) X - H \tilde{g}(X,Y)
\end{eqnarray*}
is a section of $U$, then the bundle $U$ remains totally geodesic under $\nabla$. In fact $\nabla$ is the sole preferred connection that makes $U$ and $U^{\perp}$ totally geodesic -- as adding any $\Upsilon \neq 0$ would destroy this property on at least one of these bundles.

Now we try and calculate $K$ and $K^{\perp}$ in the splitting given by $\nabla$. We know that elements of $K \cap E^{\perp}$ are of the form
\begin{eqnarray*}
\left( \begin{array}{c} 0 \\ X \\ z \end{array} \right),
\end{eqnarray*}
for some $z \in \Gamma(L^{-1})$ depending on $X$, and elements of $K \cap E^{\perp}$ are of the form
\begin{eqnarray*}
\left( \begin{array}{c} 0 \\ B \\ z' \end{array} \right).
\end{eqnarray*}
Hence, choosing $Y$ such that $Y$ and $X$ are not orthogonal,
\begin{eqnarray*}
\overrightarrow{\nabla}_Y \left( \begin{array}{c} 0 \\ X \\ z \end{array} \right) = \left( \begin{array}{c} -\mathsf{g}(Y,X) \\ \nabla_Y X -zY \\ z'' \end{array} \right)
\end{eqnarray*}
now the middle piece is a section of $U$ as well, so there exists a section
\begin{eqnarray*}
v_1 = \left( \begin{array}{c} a \\ 0 \\ z'' \end{array} \right)
\end{eqnarray*}
in $K$, with $a \neq 0$. Since $K^{\perp}$ must be orthogonal to this vector, $K^{\perp} \cap E^{\perp}$ must be of the form
\begin{eqnarray*}
\left( \begin{array}{c} 0 \\ B \\ 0 \end{array} \right),
\end{eqnarray*}
and the similar result goes for $K \cap E^{\perp}$. Consequently, as before, we have
\begin{lemm}
For any $A \in \Gamma(T)$, $\nabla_A X$ is a section of $U$.
\end{lemm}
We may, as before, choose frames $(X^j)$ and $(B^k)$ for these bundles such that the frames commute. Then
\begin{eqnarray*}
\nabla_{X^j} B^k = \nabla_{B^k} X^j = 0.
\end{eqnarray*}
This implies that the curvature tensor of $\nabla$ splits into two components, its curvature on $U$ and its curvature on $U^{\perp}$. The Ricci-tensor does the same, (see Lemma \ref{ricci:rest}), as does the rho-tensor, since $U$ and $U^{\perp}$ are orthogonal. So
\begin{eqnarray*}
\mathsf{P} = \mathsf{P}_1 + \mathsf{P}_2.
\end{eqnarray*}

We now aim to prove:
\begin{lemm}
The connection $\nabla$ is metric.
\end{lemm}
\begin{proof}
Consider the section $v_1$ in $K$, and
\begin{eqnarray*}
\overrightarrow{\nabla}_B v_1 = \left( \begin{array}{c} \nabla_B a \\ z''B + a\mathsf{P}(B) \\ \nabla_B z'' \end{array} \right).
\end{eqnarray*}
The middle term $z''B + a\mathsf{P}(B) = z''B + a\mathsf{P}_2(B)$ must be zero, showing that $\mathbf{g}^{jk}( \mathsf{P}_2)_{ij}$ is some multiple of the identity -- hence that $\mathsf{P}_2$ is a symmetric tensor. As the same is true of $\mathsf{P}_1$, $\nabla$ has symmetric rho-tensor, hence symmetric Ricci-tensor, hence preserves a volume form, hence preserves a metric $g$ in the conformal class.
\end{proof}

Defining $h_1 = g|_U$, $h_2 = g|_{U^{\perp}}$, one can, exactly as in Proposition \ref{times:manifold}, get the proof of the decomposition:
\begin{prop}
Locally, $(M,g) = (N_1, h_1) \times (N_2, h_2)$, where $N_1$ is a leaf of $U$ and $N_2$ is a leaf of $U^{\perp}$.
\end{prop}
Moreover, we've shown that $\mathsf{P}_1$ and $\mathsf{P}_2$ are multiples of $h_1$ and $h_2$ respectively; consequently $\mathsf{Ric}_1$ and $\mathsf{Ric}_2$ are as well, so both $N_1$ and $N_2$ are Einstein manifolds, with coefficients $\lambda_1$ and $\lambda_2$. We now aim to show the relation between these coefficients.

The scalar curvature $R$ of $\nabla$ is $l \lambda_1 + (n-l) \lambda_2$. Hence the rho-tensor, by Equation (\ref{confrho:formula}), is:
\begin{eqnarray*}
\mathsf{P}_1 &=& -\frac{1}{n-2} \left( \mathsf{Ric}_1 - \frac{1}{2n-2} R h_1 \right) \\
&=& -\frac{(2n-2-l) \lambda_1 + (l-n) \lambda_2 }{(n-2)(2n-2)} h_1.
\end{eqnarray*}
\begin{eqnarray*}
\mathsf{P}_2 &=& -\frac{1}{n-2} \left( \mathsf{Ric}_2 - \frac{1}{2n-2} R h_2 \right) \\
&=& -\frac{(-l) \lambda_1 + (n-2+l) \lambda_2 }{(n-2)(2n-2)} h_2.
\end{eqnarray*}

Now there is a section
\begin{eqnarray*}
v_1 = \left( \begin{array}{c} 1 \\  0 \\ f \end{array} \right)
\end{eqnarray*}
of $K$ (we may freely use $1$, as we have established that $\nabla$ is metric, hence got an isomorphism $L^{1} \cong \mathbb{R} \times M$), and a corresponding section
\begin{eqnarray*}
v_2 = \left( \begin{array}{c} 1 \\  0 \\ f' \end{array} \right)
\end{eqnarray*}
of $K^{\perp}$. Since $v_2$ is orthogonal to $v_1$, $f' = -f$. Then
\begin{eqnarray*}
\overrightarrow{\nabla}_B v_1 = \left( \begin{array}{c} 0 \\ fB - \mathsf{P}_2(B) \\ \nabla_B f \end{array} \right)
\end{eqnarray*}
as a consequence of this, we see that $f$ is a constant and
\begin{eqnarray*}
f = -\frac{(2n-2-l) \lambda_1 + (l-n) \lambda_2 }{(n-2)(2n-2)}.
\end{eqnarray*}
carrying out a similar operation on $v_2$ yields the following formula
\begin{eqnarray*}
f = \frac{(-l) \lambda_1 + (n-2+l) \lambda_2 }{(n-2)(2n-2)}.
\end{eqnarray*}
Equating these terms and re-arranging gives us the required
\begin{eqnarray*}
(n-l-1) \lambda_1 = (1-l) \lambda_2.
\end{eqnarray*}
There is, however, a rather more fundamental reason for this seemingly arbitrary equality. For:
\begin{prop}
The condition
\begin{eqnarray*}
(n-l-1) \lambda_1 = (1-l) \lambda_2.
\end{eqnarray*}
is equivalent to the rho-tensor $\mathsf{P}_{N_1}$ of $\nabla|_{N_1}$ being equal to the restriction of the rho-tensor on $M$,
\begin{eqnarray*}
\mathsf{P}_{N_1} = \mathsf{P}|_U = \mathsf{P}_1.
\end{eqnarray*}
\end{prop}
\begin{proof}
\begin{eqnarray*}
\mathsf{P}_{1} - \mathsf{P}_{N_1} &=& \left( -\frac{(2n-2-l) \lambda_1 + (l-n) \lambda_2 }{(n-2)(2n-2)} - \frac{-\lambda_1}{2(l-1)} \right) h_1 \\
&=&  \left( (n-l-1) \lambda_1 - (1-l) \lambda_2 \right) \left( \frac{(n-l)}{(l-1)(n-2)(2n-2)} \right) h_1
\end{eqnarray*}
Similarly
\begin{eqnarray*}
\mathsf{P}_{2} - \mathsf{P}_{N_2} &=& \left( -\frac{(-l) \lambda_1 + (n-2+l) \lambda_2 }{(n-2)(2n-2)} - \frac{-\lambda_2}{2(n-l-1)} \right) h_2 \\
&=&  \left( (n-l-1) \lambda_1 - (1-l) \lambda_2 \right) \left( \frac{l}{(l-1)(n-2)(2n-2)} \right) h_2
\end{eqnarray*}
Consequently, $\mathsf{P}_{1} = \mathsf{P}_{N_1}$ if and only if $\mathsf{P}_{2} = \mathsf{P}_{N_2}$, and if and only if $(n-l-1) \lambda_1 = (1-l) \lambda_2$.
\end{proof}

This is the essence of the decomposition: because of this result, $\overrightarrow{\nabla}$ operates on $\mathcal{T}_{N_1}$ along $TN_1 = U$ just as the reduced Tractor connection $\overrightarrow{\nabla}_{N_1}$ does. Now let $v_2$ be the Einstein vector in $\mathcal{T}_{N_1}$; then $\overrightarrow{\nabla}$ along $TN_2 = U^{\perp}$ will operate trivially on
\begin{eqnarray*}
K = v_2^{\perp} \cap \mathcal{T}_{N_1},
\end{eqnarray*}
since $K$ is the sum of elements of $(0,X,0)$ and $v_1$. Consequently the holonomy algebra of $\overrightarrow{\nabla}$ restricted to $K$ is $\overrightarrow{\mathfrak{hol}}_{N_1}$.

The similar result holds for $K^{\perp}$. Thus, since $K \oplus K^{\perp} = \mathcal{T}$,
\begin{eqnarray*}
\overrightarrow{\mathfrak{hol}} = \overrightarrow{\mathfrak{hol}}_{N_1} \oplus \overrightarrow{\mathfrak{hol}}_{N_2}.
\end{eqnarray*}

To reverse this decomposition, define $(M,g)$ as $(N_1,h_1) \times (N_2, h_2)$ with $N_1$ and $N_2$ Einstein with Einstein coefficients related as above. Then the overall Tractor connection $\overrightarrow{\nabla}$ will be generated by $\overrightarrow{\nabla}_{N_1}$ and $\overrightarrow{\nabla}_{N_2}$ as above. Then let $v_2$ be the Einstein vector of $\mathcal{T}_{N_1}$. Then the bundle
\begin{eqnarray*}
K = v_2^{\perp} \cap \mathcal{T}_{N_1},
\end{eqnarray*}
is preserved by $\overrightarrow{\nabla}$ as is its orthogonal complement
\begin{eqnarray*}
K = v_1^{\perp} \cap \mathcal{T}_{N_2},
\end{eqnarray*}
where $v_1$ is the Einstein vector of $\mathcal{T}_{N_2}$. Note that $v_1 \in \Gamma(K)$ and $v_2 \in \Gamma(K^{\perp})$, which explains the somewhat odd numbering of them.

\begin{exa} To illustrate these proofs, we can see that $S^4 \times \mathbb{R}^4$ does not have any holonomy-conserved sub-bundles in its tractor connection (in fact it has full holonomy), while $S^4 \times \mathbb{H}^4$ is conformally flat, for $\mathbb{H}^4$ the hyperbolic 4-space. \end{exa}

\begin{rem} Some old results of H.~W.~Brinkmann \cite{RSCES}, \cite{ESMC} can be proved directly using this decomposition theorem. For instance, the fact that any 4-manifold with two distinct Einstein structures in the conformal class is conformally flat (a direct consequence of the flatness of any M\"obius structure with reduced holonomy, see next section). In our setting, the preserved sub-bundle spanned by the two Einstein vectors decomposes the manifold into a direct product of 3- and 1-dimensional Einstein spaces. But both these spaces are conformally flat, so our original manifold has trivial holonomy; in other words, it is conformally flat. \end{rem}

\begin{rem} Analogously to the previous remark, we can see that if not conformally flat, a five dimensional manifold can have up to two linearly independent Einstein structures, a six dimensional manifold can have three, an $n$ dimensional manifold $n-3$.\end{rem}

\section{Einstein Spaces: Metric Cones}

In this section we will give a full classification of the possible Tractor holonomies of the non Ricci-flat Einstein spaces, using to this effect the construction of a metric cone, whose Levi-Civita holonomy corresponds to the Tractor holonomy of the original manifold.

\begin{rem} As we mentioned in the introduction, this metric cone construction is related to the Ambient Metric construction of \cite{CI} and \cite{STCAMC}, for conformally Einstein manifolds. The actual relation is slightly subtle. This also provides a direct proof of a result of \cite{AOT}, namely that the Ambient Metric construction always exists if the manifold is conformally Einstein. \end{rem}

\begin{defi} A conformal manifold $M$ is said to be \emph{indecomposable} if they cannot be decomposed into Einstein spaces as in the previous section. In other words, $\overrightarrow{\nabla}$ may preserve a single line bundle (and its orthogonal complement), but nothing else. \end{defi}

\begin{nrem} \label{irr:hol} In the non Ricci-flat Einstein case, indecomposable implies that the tractor holonomy acts irreducibly on $\mathbb{R}^{n+1}$ or $\mathbb{R}^{n,1}$ (since the only preserved line bundle is positive or negative definite). \end{nrem}

As all Einstein manifolds of dimension 3 are conformally flat, we shall assume $n > 3$.

Let $(M,g)$ be an Einstein manifold, $\mathsf{Ric} = \lambda g$, $\lambda \neq 0$.

\begin{theo}[Einstein Classification] \label{ein:cla} The Tractor holonomy of $M^n$ is one of the following, $n \geq 4$:
\begin{itemize} \item[-] $SO(n,1)$, \item[-] $SO(n+1)$, \item[-] $SU(m)$ for $2m=n+1$, \item[-] $Sp(m)$ for $4m=n+1$, \item[-] $G_2$ for $n=6$, \item[-] $Spin(7)$ for $n=7$. \end{itemize}
Moreover, all these holonomy groups actually occur.
\end{theo}

\begin{rem} It is interesting to note that there is only a single holonomy possible for an indecomposable Einstein manifold with negative constant. \end{rem}

The remainder of this chapter will be dedicated to proving Theorem \ref{ein:cla}.

\begin{defi}
Given an Einstein manifold $(M, g)$, we define the metric cone on $M$ as $(N = \mathbb{R}^+ \times M, h)$ with
\begin{eqnarray*} h = \frac{1}{\mu} dt^2 + t^2 g,\end{eqnarray*}
and $\mu = \frac{\lambda}{n-1}$.
\end{defi}

Note that $h$ is of definite signature if and only if $M$ has positive Einstein constant. In the negative case, we call $(N, h)$ a Lorentzian cone.

Then defining $\nabla$ as the Levi-Civita connection of $N$, and remembering the formula:
\begin{eqnarray*} 2\langle\nabla_{X} Y, Z\rangle & = & X . \langle Y, Z\rangle + Y . \langle X, Z\rangle - Z . \langle X, Y\rangle \\ & & + \langle[X, Y], Z\rangle + \langle [Z, X], Y\rangle + \langle [Z, Y], X\rangle, \end{eqnarray*}
we can calculate the following equalities. For $T = \frac{\partial}{\partial t}$, and $X_i$ a local basis of vector fields of $M$, extended trivially to $N$:
\begin{eqnarray*} \nabla_T T &=& 0, \\ \nabla_{X_i} T &=& \frac{1}{t} X_i, \\ \nabla_T X_i &=& \frac{1}{t} X_i, \\ \nabla_{X_i} X_j &=& \widetilde{\nabla}_{X_i} X_j - t \mu g(X_i, X_j) T,  \end{eqnarray*} 
with $\widetilde{\nabla}$ the Levi-Civita connection of $g$.

Given a path $\tau$ in $\{1\} \times M$, with tangent vector field $Z$, let $Y$ be the parallel transport of a vector along the path, thus $\nabla_Z Y = 0$. Split $Y$  as $Y^{\perp} + aT$, with $Y^{\perp} \in \Gamma (T)$. Then we get the following result: 

\begin{lemm} Extend $Z$ and $Y$ in the $T$ direction, with $Z(t,x) = Z(x)$, and 
\begin{eqnarray*} Y(t,x) = \frac{1}{t} Y^{\perp}(x) + a(x)T. \end{eqnarray*}
Then $\nabla_T Y = 0$ and $\nabla_Z Y = 0$ on $\tau \times \mathbb{R}^+$. \end{lemm}

\begin{proof}

The function $a$ is independent of $t$, so $T(a) = 0$. Hence $\nabla_T (aT)=0$. Furthermore,
\begin{eqnarray*}
\nabla_T \left( \frac{1}{t} Y^{\perp} \right) &=& T \left(\frac{1}{t} \right) Y^{\perp} + \frac{1}{t^2} Y^{\perp} \\ &=& (-\frac{1}{t^2} + \frac{1}{t^2}) Y^{\perp} \\ &=& 0,
\end{eqnarray*}
so $\nabla_T Y = 0$.

We can expand out the original equation $\nabla_Z Y = 0$ at $t=1$, giving:
\begin{eqnarray*} 0 = \widetilde{\nabla}_{Z} Y^{\perp} - \mu g(Z, Y^{\perp}) T + aZ + Z(a)T.\end{eqnarray*}

By linearity, this is equivalent to the two equations $0 = \widetilde{\nabla}_{Z} Y^{\perp} + aZ$ and $0=(-\mu g(Z, Y^{\perp})+ Z(a))T$.

Then similarly expanding $\nabla_Z Y$ for varying $t$:
\begin{eqnarray*} \nabla_Z Y &=& \frac{1}{t} \widetilde{\nabla}_{Z} Y^{\perp} - \frac{t}{t} \mu g(Z, Y^{\perp}) T + \frac{a}{t}Z + Z(a)T \\ &=& \frac{1}{t}(\widetilde{\nabla}_{Z}Y^{\perp} + aZ) + (-\mu g(Z, Y^{\perp}) + Z(a))T \\ &=& 0.\end{eqnarray*}
\end{proof}

The previous result shows that when we're computing the holonomy of $\nabla$, we only need to consider paths in $\{1\} \times M \cong M$.

We can now turn to the tractor connection $\overrightarrow{\nabla}$ on $\mathcal{T}$, for the conformal structure $\mathbf{g} \simeq [g]$. Using the splitting given by the metric $g$, we can see the formal similarities with $\nabla$ at $t = 1$.

As $g$ is Einstein, with coefficient $\lambda$, then $\overrightarrow{\nabla}  \left( \frac{n-1}{\lambda}, 0 , -\frac{1}{2} \right) = 0$. Furthermore, for $R = \left( \frac{n-1}{\lambda}, 0, \frac{1}{2} \right)$, then:

\begin{eqnarray*} \overrightarrow{\nabla}_{X_i} R = \left( \begin{array}{c} 0 \\ X_i \\ 0 \end{array} \right) \ \ \mathrm{and} \  \ \overrightarrow{\nabla}_{X_i} \left( \begin{array}{c} 0 \\ X_j \\ 0 \end{array} \right) = \left( \begin{array}{c} 0 \\ \widetilde{\nabla}_{X_i} X_j \\ 0 \end{array} \right) - \mu g(X_i, X_j) R. \end{eqnarray*}

Hence under the formal identification of $R$ with $T$ and $(0, X_i, 0)$ with $X_i$, we get $\nabla_Z \cong \overrightarrow{\nabla}_Z$ for $Z \in \Gamma(T)$ at $t=1$. Then by the previous lemma and its implication for the holonomy of $\nabla$, the next theorem is proved:

\begin{theo} The holonomy groups of $(\mathcal{T}, \overrightarrow{\nabla}, (M, \mathbf{g}))$ and $(TN, \nabla, (N,h))$ are isomorphic. \end{theo}

Hence the holonomy of $\overrightarrow{\nabla}$ is metric, and irreducible by Remark \ref{irr:hol}, and must be one of those classified by Merkulov and Schwachh\"ofer in \cite{CIH}. In the negative Einstein case, a look at the table shows that the only possible holonomy is the full $SO(n,1)$ group itself. For the positive Einstein, we need the following result:

\begin{prop} The metric cone $(N, h)$ is Ricci-flat.\end{prop}

\begin{proof}
By the definition of $h$ and the corresponding $\nabla$, the curvature $R = R^\nabla$ is
\begin{eqnarray*} R_{T,-} &=& 0, \\ R_{X_i X_j} X_k &=& \widetilde{R}_{X_i X_j} X_k - \mu g(X_j, X_k) X_i + \mu g(X_i, X_k) X_j, \end{eqnarray*}
with $\widetilde{R}$ the curvature of $\widetilde{\nabla}$.

Then taking traces,
\begin{eqnarray*} \mathsf{Ric} (T,-) &=& 0, \\ \mathsf{Ric}(X_j, X_k) &=& \widetilde{\mathsf{Ric}}(X_j,X_k) + (1-n) \mu g(X_j, X_k) \\ &=& \lambda g(X_j, X_k) - \lambda g(X_j, X_k) \\ &=& 0. \end{eqnarray*}
\end{proof}

So the possible holonomies reduce to those corresponding to metrics which are Ricci-flat, namely $SO(n+1), SU(m), Sp(m), G_2$ and $Spin(7)$.

The $SO(n+1)$ case is generic. The $SU(m)$ holonomy on the cone corresponds to Sasaki-Einstein manifolds, the $Sp(m)$ to 3-Sasakian ones, $G_2$ and $Spin(7)$ to weak holonomy manifolds \cite{RKSH}; all of which can be realised on compact manifolds.

It is immediate that a metric cone on any one dimensional space is flat. We now aim to show that the metric cone on a two dimensional Einstein space of constant scalar curvature is also flat.
\begin{prop} \label{two:flat}
Any Tractor connection in two dimensions with a preserved Tractor $u$ is flat.
\end{prop}
\begin{proof}
In this case, we have
\begin{eqnarray*}
\mathsf{P} = \mu Id,
\end{eqnarray*}
With $\mu$ a constant.
However, \cite{Mob}, the only curvature element of a Tractor/M\"obius connection in two dimensions is the Cotton-York tensor -- which must vanish entirely, as $\nabla \mathsf{P} = 0$, making the connection flat.
\end{proof}

\begin{rem} In \cite{LTT} R.J.~Baston presents a local twistor theory, which, in the case of conformal manifolds, is just given by the spin representation of $\mathcal{G}_0$ and the extension of $\overrightarrow{\nabla}$ to this new context. A parallel section of this bundle is equivalent with the existence of a spinor $\psi$ solving the \emph{twistor equation} for all vector field $X$:
\begin{eqnarray*} \nabla_X \psi + \frac{1}{n} X. D \psi = 0, \end{eqnarray*}
with $D$ the Dirac operator. Paper \cite{TERM} by Katharina Habermann analyses solutions to this twistor equation; she shows that these imply that the manifold is conformally Einstein, of non-negative scalar curvature.

So the Tractor holonomy groups $G_2$ and $Spin(7)$ actually correspond to the existence of twistor-spinors on the manifold.\end{rem}

\begin{rem} The concept of a twistor-spinor is a generalisation of that of a Killing spinor. A Killing spinor is a spinor $\psi$ solving the equation
\begin{eqnarray*}
\nabla_X \psi = \lambda X . \psi
\end{eqnarray*}
for all vector fields $X$ and some constant $\lambda$. In \cite{RKSH}, C.~B\"ar showed that having a Killing spinor is equivalent with having a parallel spinor on the metric cone. So the cases of weak holonomy SU(3) and nearly K\"ahlerian structures are covered by the Tractor connection; in fact in his paper \cite{MEH} constructing manifolds of exceptional holonomy, R.L.~Bryant produces manifolds of holonomy $G_2$ and $Spin(7)$ as metric cones on $SU(3)/T^2$ and $SO(5)/SO(3)$ respectively. Thus all the holonomies listed actually occur. \end{rem}

We can now turn to the Ricci-flat case, which is actually simpler than the general Einstein case, but with an added subtlety.

\section{Ricci-Flat Spaces}

Let $(M^n,g)$ be a Ricci-flat space of indecomposable tractor holonomy. As $M^n$ is Ricci-flat, its Tractor holonomy is contained within $SO(n) \rtimes \mathbb{R}^n$. Fix a point $b \in M$ for calculating the holonomy groups, and let $H$ be the metric holonomy of $M$, $D$ its Tractor holonomy, $\mathfrak{h}$, $\mathfrak{l}$ their Lie algebras.

Then:

\begin{lemm} $H \subset D$, or, equivalently, $\mathfrak{h} \subset \mathfrak{l}$. \end{lemm}

\begin{proof}
Let $Y$ be the parallel transport of a vector along a path $\tau$ with tangent field $X$; in other words $\nabla_X Y =0$, for $\nabla$ the metric connection on $M$. Then
\begin{eqnarray*} \overrightarrow{\nabla}_X \left( \begin{array}{c} x \\ Y \\ 0 \end{array} \right) = \left( \begin{array}{c} \nabla_X x - g(X,Y) \\ \nabla_X Y = 0 \\ 0 \end{array} \right) , \end{eqnarray*}
which is zero for $x = \int_{\tau} g(X,Y)$, proving that every metric holonomy element is a tractor holonomy element. This argument also works in reverse, showing that $\pi (\mathfrak{l}) = \mathfrak{h}$, where $\pi$ is the projection of $\mathfrak{co}(n)_b \oplus T_b$ onto its first component. \end{proof}

This demonstrates that $\mathfrak{l} \subset \mathfrak{b} = \mathfrak{h} \oplus T_b$. But first:

\begin{lemm} The representation of $\mathfrak{h}$ on $T_b$ is irreducible. \end{lemm}

\begin{proof}

If a bundle $S \subset T$ is preserved by $\nabla$, then the bundle $\mathcal{E}[-1] \oplus S$ is preserved by $\overrightarrow{\nabla}$. Thus, since we assume our Tractor holonomy to be indecomposable, then $\mathfrak{h}$ must act irreducibly on $T_b$. \end{proof}

Then since the Lie bracket on $T_b$ is trivial, the adjoint representation of $\mathfrak{h} \subset \mathfrak{l} \subset \mathfrak{b}$ on the second component of $\mathfrak{b}$ is the usual, irreducible one. Accordingly this adjoint representation splits $\mathfrak{b}$ into two irreducible representations, isomorphic to $\mathfrak{h}$ and $T_b$.

As a consequence, $\mathfrak{l} \cong \mathfrak{h}$ or $\mathfrak{l} \cong \mathfrak{b}$. We now claim that

\begin{lemm} $\mathfrak{l} \cong \mathfrak{b}$. \end{lemm}

\begin{proof}

Reasoning by contradiction, we assume that $\mathfrak{l} \cong \mathfrak{h}$, and go on to show that this violates our indecomposability assumption.

Express $\mathfrak{b}$ as $\mathfrak{h}_0 \oplus \mathbb{R}^n_0$, the sum of the irreducible representations of $\mathfrak{h}$. Then, as $\mathfrak{h}_0 \cong \mathfrak{h}$ acts irreducibly on $T_b$, there is, at $b$, a new splitting of $\mathcal{T}$ corresponding to the splitting $\mathbb{R}^{n*} \oplus \mathfrak{h}_0 \oplus \mathbb{R}^{n}_0$. This splitting is
\begin{eqnarray*} \mathcal{E}[-1]_0 \oplus T[-1]_0 \oplus \mathcal{E}[-1] \end{eqnarray*}

Then $\mathfrak{h} \subset \mathfrak{l}$ preserves the new vectors $(1, 0, 0)$ and $(0, 0, 1)$. This shows that $\mathfrak{l}$ preserves a rank two sub-bundle, contradicting indecomposability. \end{proof}

Putting this together, we can now claim the following theorem:

\begin{theo} The possible indecomposable Tractor holonomy groups for the conformal manifold $(M, \mathbf{g})$, conformally Ricci-flat, are:
\begin{itemize}
\item[-] $SO(n) \rtimes \mathbb{R}^n, n \geq 4$ (generic),
\item[-] $SU(m) \rtimes \mathbb{R}^{2m}, m \geq 2$ (Calabi-Yau),
\item[-] $Sp(m) \rtimes \mathbb{R}^{4m}, m \geq 1$ (Hyper-K\"ahler),
\item[-] $G_2 \rtimes \mathbb{R}^7$ (see \cite{MEH}),
\item[-] $Spin(7) \rtimes \mathbb{R}^8$ (see \cite{MEH}),
\end{itemize}
and all of these groups do occur.
\end{theo}

\begin{rem}
This result offers an alternative proof for the theorems in Mario Listing's paper \cite{CEPND}, in the special case of conformally Ricci-flat manifolds.
\end{rem}
\begin{rem}
The metric cones constructed in the previous chapter are Ricci-flat (pseudo-)Riemannian manifolds. They are not, however, indecomposable; in fact, their Tracor holonomy is equal to their metric holonomy. This property caracherises metric cones.
\end{rem}

\section{Addendum: Symmetric Spaces}

A symmetric space $(S, g)$ is a manifold such that $\nabla^{g} R^g = 0$ for $R^g$ the full curvature tensor. It is quite easily to show, using the infinitesimal holonomy developed by S.~Kobayashi and K.~Nomizu \cite{KobNom}, that any indecomposable conformal manifold that is conformal to a symmetric space has the maximal holonomy in its category - $SO(n+1,1)$ if the symmetric space is not Einstein, and $SO(n+1)$ or $SO(n,1)$ if it is (Ricci-flat symmetric spaces are flat).

\begin{exa} These results give an independent proof to the results of F.~Leitner \cite{CHBM}, that the conformal holonomy of $SO(4)$, locally isomorphic to $S^3 \times S^3$, is $SO(7)$. The group $SO(4)$ is a positive Einstein symmetric space, not conformally flat (consider the fate of the tractor vector $(1,0,0)$ under parallel translation), so the result follows. \end{exa}

\begin{exa} Note that the same argument shows that the manifold $S^n(a) \times S^n(b)$ where $a \neq b$ are the radii of the spheres, has full holonomy $SO(2n+1,1)$. \end{exa}

Two very similar results also are implied:
\begin{exa} If a manifold $(M, \mathbf{g})$ is conformal to an Einstein symmetric space, then it cannot be conformal to any other Einstein space, or any other symmetric Space, unless it is conformally flat. \end{exa}
And:
\begin{exa} If $(M, \mathbf{g})$ is conformal to a symmetric space in two different ways, then its tractor holonomy is full or null. \end{exa}

\section{Indefinite Signature}
Most of the results of this paper extend to the general pseudo-Riemannian case, see \cite{methesis}. The cone construction and the Ricci-flat results still apply, as do the results about umbilicity. The decomposition theorem, however, requires an extra condition: that $K \cap K^{\perp}$ be of rank one or zero (a condition that is automaticaly true in the definite signature case). This is equivalent with requiring that $U \cap U^{\perp} = 0$; without it, the decomposition can't proceed.

And, unlike the result for $\mathfrak{so}(n+1,1)$ proved in \cite{GHSHS}, there are non-trivial subalgebras of $\mathfrak{so}(p+1,q+1)$ acting irreducibly on $\mathbb{R}^{(p+1,q+1)}$. Thus we have many other candidate algebras to deal with.

\end{document}